\documentclass[12pt, reqno]{amsart}
\usepackage{latexsym}
\usepackage{mathtools}
\usepackage{amsmath}
\usepackage{amsthm}
\usepackage{mathrsfs}
\usepackage{lmodern}
\usepackage[T1]{fontenc}
\usepackage{textcomp}
\usepackage[utf8]{inputenc}
\usepackage{color}
\usepackage{amssymb}
\usepackage{enumerate}
\usepackage[abbrev]{amsrefs}
\usepackage{indentfirst}
\usepackage{graphicx}
\usepackage{bmpsize}
\usepackage{hyperref}
\usepackage{enumitem}
\usepackage{enumerate}
\usepackage{url}
\usepackage{autobreak}
\usepackage[justification=justified]{caption}
\usepackage{bm}
\usepackage{here}
\usepackage[labelformat=empty,subrefformat=parens]{subcaption}
\renewcommand{\thefootnote}{} 

\theoremstyle{plain} 
\newtheorem{theorem}{\indent\sc Theorem}[section]
\newtheorem{lemma}[theorem]{\indent\sc Lemma}

\theoremstyle{definition} 
\newtheorem{definition}[theorem]{\indent\sc Definition}
\newtheorem{remark}[theorem]{\indent\sc Remark}
\newtheorem{example}[theorem]{\indent\sc Example}

\newtheorem{convention}[theorem]{\indent\sc Convention}
\newtheorem{question}[theorem]{\indent\sc Question}
\makeatletter
  
  \@addtoreset{equation}{section}
\makeatother

\newcommand{\cF}{\mathcal{F}}
\newcommand{\cL}{\mathcal{L}}
\newcommand{\cT}{\mathcal{T}}
\newcommand{\cB}{\mathcal{B}}

\newcommand\cinput[2]{\lower#1pt\hbox{\input{#2}}}

\makeatletter
\@namedef{subjclassname@2020}{\textup{2020} Mathematics Subject Classification}
\makeatother

\begin{document}

\title[The $3$-colorable subgroup $\mathcal{F}$ and tricolorability of links]{THE $3$-COLORABLE SUBGROUP OF THOMPSON'S GROUP AND TRICOLORABILITY OF LINKS}
\author{Yuya Kodama and Akihiro Takano}
\date{}
\renewcommand{\thefootnote}{\arabic{footnote}}  
\setcounter{footnote}{0} 
\thanks{The first author was supported by JST, the establishment of university
fellowships towards the creation of science technology innovation, Grant Number JPMJFS2139.}
\begin{abstract}
Starting from the work by Jones on representations of Thompson's group $F$, subgroups of $F$ with interesting properties have been defined and studied. 
One of these subgroups is called the $3$-colorable subgroup $\mathcal{F}$, which consists of elements whose ``regions'' given by their tree diagrams are $3$-colorable.
On the other hand, in his work on representations, Jones also gave a method to construct knots and links from elements of $F$. 
Therefore it is a natural question to explore a relationship between elements in $\mathcal{F}$ and $3$-colorable links in the sense of knot theory. 
In this paper, we show that all elements in $\mathcal{F}$ give 3-colorable links. 
\end{abstract}
\maketitle
\section{Introduction}
Thompson's groups were defined by Richard Thompson from the motivation of logic. 
Later these groups turned out to be related to various areas and have been studied from many viewpoints, such as algebraic and analytic properties, subgroups, quasi-isometric invariants, and so on. 

Vaughan Jones \cite{jones2017thompson} began the study of several representations of Thompson groups $F$ and $T$ motivated by the creation of algebraic quantum field theories on the circle. 
In this process, two interesting research projects of $F$ arose besides the study of representations: 
\begin{enumerate}
\item study of new subgroups of $F$, and
\item study of relationships between elements of $F$ and knots and links obtained from them. 
\end{enumerate}
In (1), principal groups are the group called oriented subgroup $\overrightarrow{F}$ \cite{jones2017thompson} and the group called $3$-colorable subgroup $\cF$ \cite{jones2018nogo}, where $\cF$ is the group we will focus on in this paper. 
These groups are each isomorphic to a certain Brown--Thompson group \cite{golan2017jones, ren2018skein}. 
Also, each of these groups is a stabilizer of a subset of the closed interval $[0, 1]$ under the natural action when each element is regarded as a piecewise linear map on $[0, 1]$ \cite{golan2017jones, aiello2021maximal}. 
This implies they are closed in $F$. 
In contrast, it is known that there exists a maximal subgroup of infinite index in $F$ isomorphic to $\overrightarrow{F}$ such that it does not fix any point in the open interval $(0, 1)$ \cite{golan2017subgroups}. 
For (2), it was already shown that any (oriented) link is obtained from an element of $F$ (resp.~ $\overrightarrow{F}$) \cite{aiello2020alexander, jones2017thompson}, which is called Alexander's theorem.
Also it is known that oriented links obtained from ``positive words'' in $\overrightarrow{F}$ are positive \cite{aiello2021positive}, and links obtained from ``positive words'' in $F$ are arborescent \cite{aiello2022arborescence}. 

On the other hand, to the best of the authors' knowledge, there has been no study of $\cF$ from the viewpoint of knot theory. 
Indeed, in Aiello's survey \cite{aiello2022thompson}, he made the same claim and asked the following question: 
\begin{question}[{\cite[Question 4]{aiello2022thompson}}]
Do the elements of $\cF$ produce all unoriented knots and links?
\end{question}

We will answer this question negatively by showing the following theorem. 
\begin{theorem}\label{main}
All knots and links obtained from non-trivial elements of $\cF$ are $3$-colorable.
\end{theorem}
Since there exist non $3$-colorable knots and links such as the figure-eight knot, and Whitehead link, we get the desired result. 

\section{Preliminaries} 
\subsection{Thompson's group $F$ and the $3$-colorable subgroup $\cF$} \label{thompson}
We summarize the definitions of Thompson's group $F$ and its subgroup $\cF$ called $3$-colorable subgroup. 
It is known that there exist several (equivalent) definitions of Thompson's group $F$. 
We use pairs of binary trees in this paper.
See \cite[\S 2]{cannon1996intro} for details. 

Let $\cT$ be a set consisting of pairs of binary trees with the same number of leaves. 
We define an equivalence relation on $\cT$ as follows: 
let $(T_+, T_-) \in \cT$ with $n$ leaves. 
We label the leaves of the trees with $1, \dots, n$ from left to right. 
Assume that there exists $i \in \{1, \dots, n\}$ such that two leaves labeled $i$ and $i+1$ have a common parent in $T_+$ and $T_-$, respectively. 
Then we can remove a minimal binary tree containing two leaves labeled $i$, $i+1$, and the common parent from $T_+$, and obtain a binary tree $T_+^\prime$. 
Similarly, we obtain $T_-^\prime$. 
Therefore we get an element $(T_+^\prime, T_-^\prime)$ in $\cT$ from $(T_+, T_-)$. 
We call this operation and its inverse operation \textbf{reduction} and \textbf{insertion of carets}, respectively. 
Then we define an equivalence relation $\sim$ as the one generated by these operations. 
We say $(T_+, T_-)$ is \textbf{reduced} if we can not reduce any carets. 
It is known that there exists a unique reduced representative for each equivalence class \cite[\S 2]{cannon1996intro}. 

\begin{figure}[tbp]
\begin{center}
\includegraphics[height=200pt]{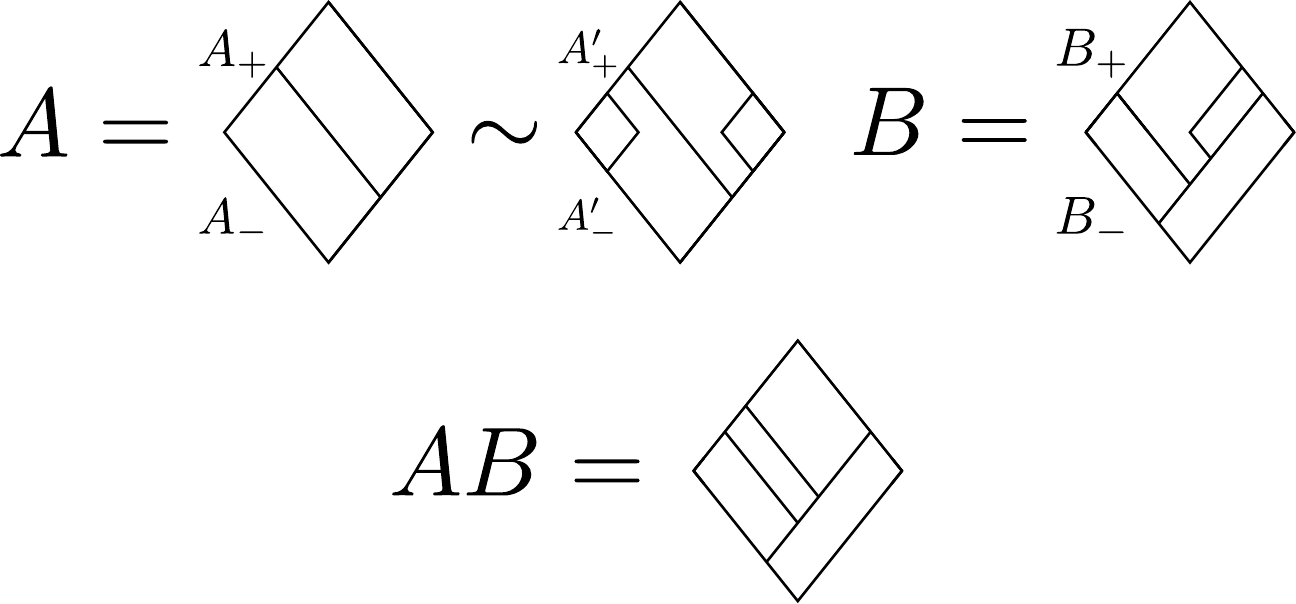}
\end{center}
\caption{Example of the product. }
\label{product_example}
\end{figure}

We call \textbf{Thompson's group $F$} the group $(\cT/{\sim}, \times)$ where $\times$ is the following product. 
Let $A=(A_+, A_-)$ and $B=(B_+, B_-)$ be in $\cT$. 
By insertions, we can obtain $A^\prime=(A_+^\prime, A_-^\prime)$ and $B^\prime=(B_+^\prime, B_-^\prime)$ such that $A \sim A^\prime$, $B \sim B^\prime$, and $A_-^\prime=B_+^\prime$ hold. 
Then for the equivalence classes of $A$ and $B$, we define their product as an equivalence class of $(A_+^\prime, B_-^\prime)$. 
Figure \ref{product_example} is an example of this product. 
\begin{convention}
As in Figure \ref{product_example}, we always assume that for each $(T_+, T_-)$ in $\cT$, $T_+$ is descending tree with the root on top, $T_-$ is ascending with root on lowest, and that leaves with the same label are attached. 
We call this diagram \textbf{tree diagram}. 
\end{convention}

In the rest of this section, we discuss the $3$-colorable subgroup $\cF$.
\begin{definition}
For a tree diagram, put it in $\mathbb{R}^2$ so that the top root is on $(0, 1) \in \mathbb{R}^2$ and the lowest root is on $(0, -1) \in \mathbb{R}^2$.
Also we add the short edges connecting the top root $(0,1)$ and $(0,2)$ (resp. the bottom root $(0,-1)$ and $(0,-2)$).
Then we can obtain a partition of $\mathbb{R} \times [-2, 2]$ by using the edges of the tree diagram. 
This partition is called the \textbf{strip} of this tree diagram. 
\end{definition}

\begin{figure}[tbp]
\begin{center}
\includegraphics[height=165pt]{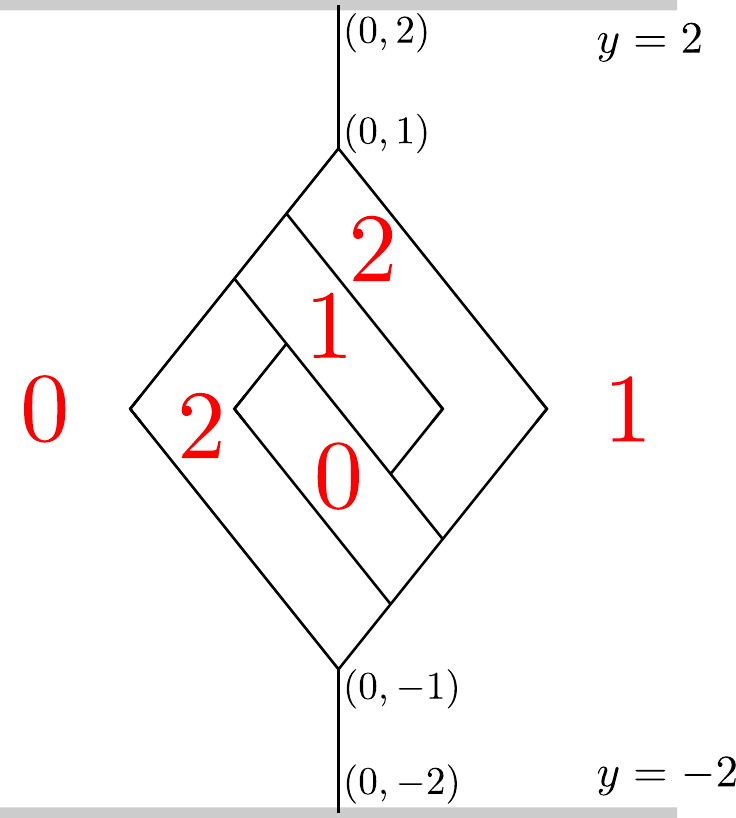}
\end{center}
\caption{A tree diagram in $\mathbb{R} \times [-2, 2]$ with regions colored by red letters. }
\label{3_colorable_element}
\end{figure}

We define $\cF$ by the following $3$-colorability of strips. 
\begin{definition}
A tree diagram is \textbf{$3$-strip-colorable} if we can assign the colors $\{0, 1, 2\}$ to each partitioned region of the strip such that if two regions have a common edge, then they have conflicting colors. 
By convention, if a tree diagram is $3$-strip-colorable, then assign $0$ to the left side and $1$ to the right side of the two infinite regions. 
See Figure \ref{3_colorable_element} for an example of a $3$-strip-colorable tree diagram. 
\end{definition}

\begin{figure}[tbp]
\begin{center}
\includegraphics[height=20pt]{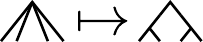}
\end{center}
\caption{Replacement rule. }
\label{4adic_to_2adic}
\end{figure}

\begin{figure}[tbp]
\begin{center}
\includegraphics[height=150pt]{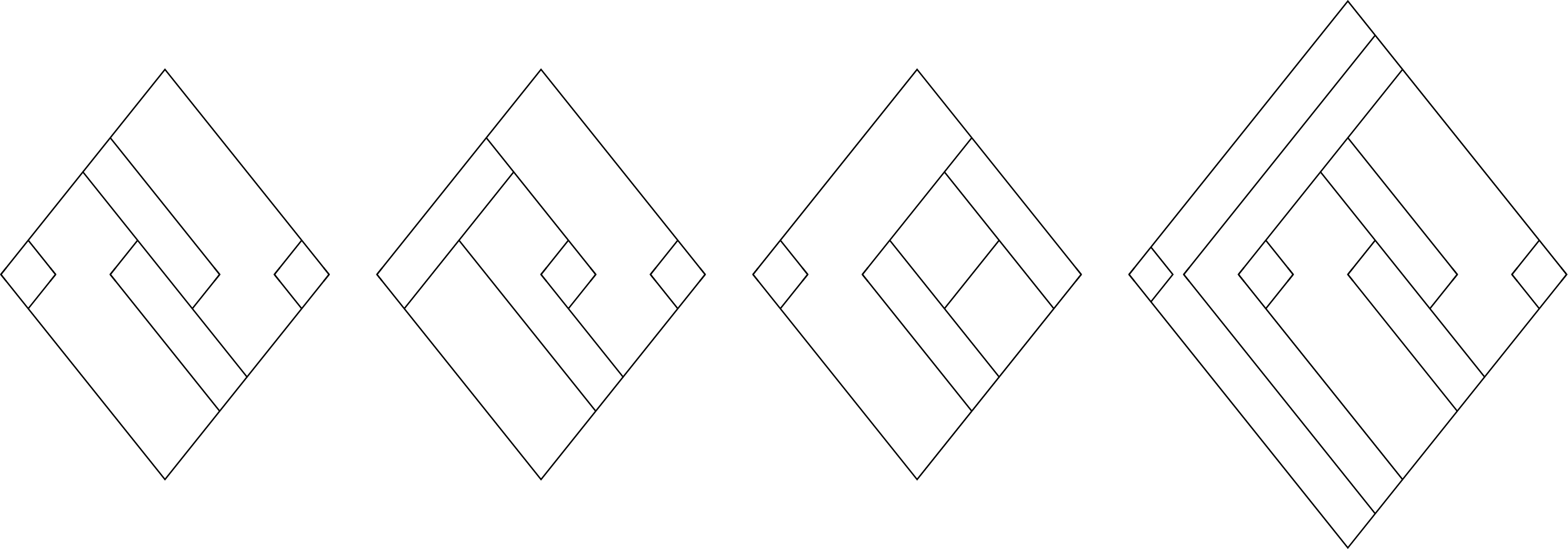}
\end{center}
\caption{The (non-reduced) generators $w_0$, $w_1$, $w_2$, and $w_3$. }
\label{generator_colF}
\end{figure}

By the above convention, if a tree diagram is $3$-strip-colorable, then its coloring is unique. 
Also, note that if a tree diagram is $3$-strip-colorable, then any tree diagram that is equivalent to it is also $3$-strip-colorable. 
Now we define the group $\cF$. 
\begin{definition}
We define the $3$-colorable subgroup $\cF$ as a subgroup of Thompson's group $F$ consisting of equivalence classes of $3$-strip-colorable tree diagrams. 
\end{definition}
This group is obtained from an embedding of Brown--Thompson group $F(4)$ to $F$. 
This embedding is defined as ``caret replacements'', where the replacement rule is illustrated in Figure \ref{4adic_to_2adic}. 
See \cite[Sections 1 and 3]{burillo2001metrics} for the precious definition of $F(4)$ and this embedding. 
\begin{theorem}[{\cite[Theorem 1.3]{ren2018skein}}] \label{Thm_cF_Brown-Thompson}
The group $\cF$ is isomorphic to Brown--Thompson group $F(4)$. 
Especially, $\cF$ is generated by $w_0$, $w_1$, $w_2$, and $w_3$ depicted in Figure $\ref{generator_colF}$. 
\end{theorem}
\subsection{Jones' construction}
We recall the method of constructing links from elements of $F$ by referring to \cite{jones2019thompson}.
Let $(T_+, T_-)$ be a reduced tree diagram.

\noindent
\textbf{Step 1: Construct the plane graph $\cB(T_+, T_-)$.}

\begin{figure}[tbp]
\begin{center}
\includegraphics[height=135pt]{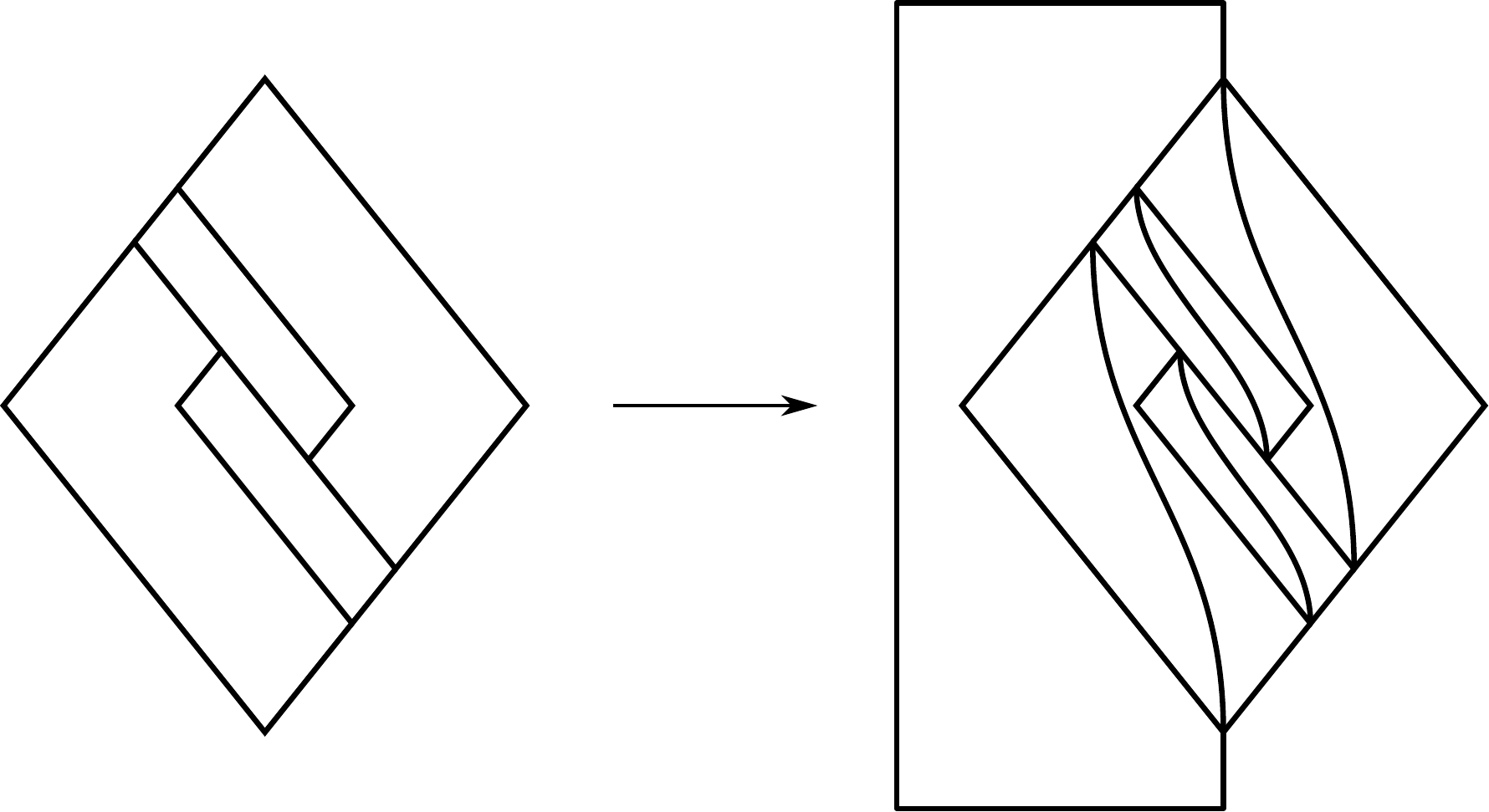}
\end{center}
\caption{The plane graph $\cB(T_+, T_-)$ associated with $(T_+, T_-)$.}
\label{step1}
\end{figure}

The plane graph $\cB(T_+, T_-)$ associated with a tree diagram $(T_+, T_-)$ is defined as follows:
each region, including the unbounded one, of $(T_+, T_-)$ has a unique vertex of $T_+$ (resp.~ $T_-$) such that its bifurcating edges are part of edges of the region.
Then we join such two vertices with an edge in the region; see Figure \ref{step1}.

\noindent
\textbf{Step 2: Construct the link diagram $\cL(T_+, T_-)$.}

\begin{figure}[tbp]
\begin{center}
\includegraphics[height=135pt]{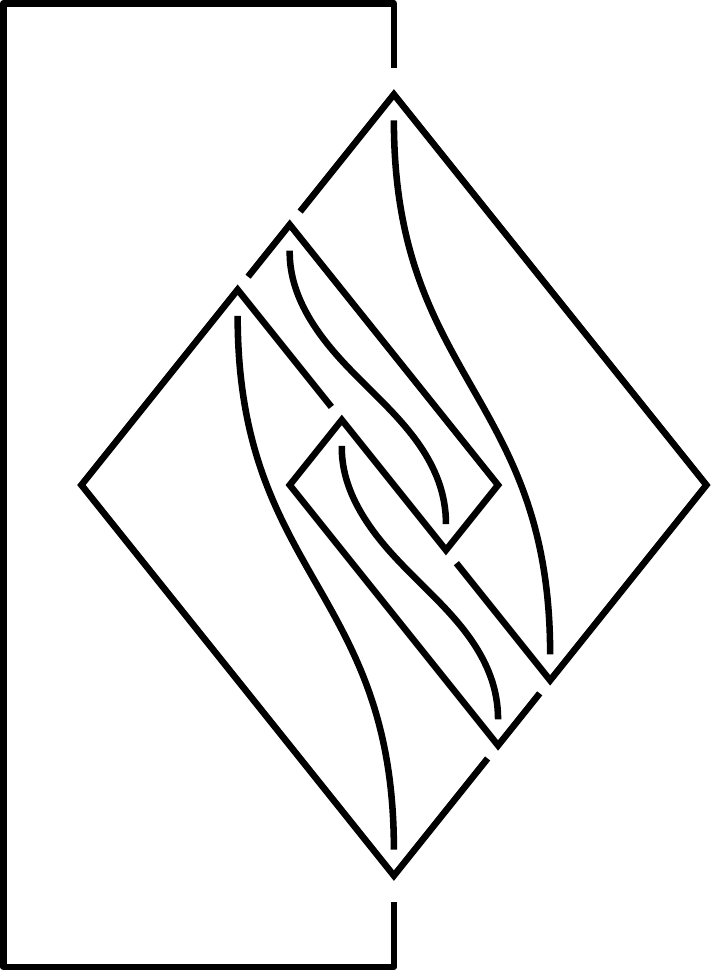}
\end{center}
\caption{The link diagram $\cL(T_+, T_-)$.}
\label{step2}
\end{figure}

Since all vertices of the graph $\cB(T_+, T_-)$ are 4-valent, we are able to regard it as a link projection.
Then we turn them into crossings with the rules \cinput{5}{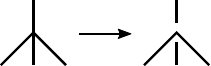} or \cinput{5}{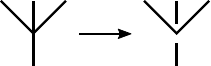}, and obtain the link diagram $\cL(T_+, T_-)$; see Figure \ref{step2}.

\begin{figure}[tbp]
\begin{center}
\includegraphics[height=50pt]{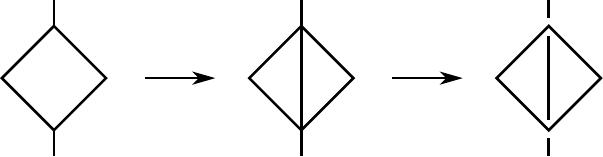}
\end{center}
\caption{A trivial link component by a caret.}
\label{caret_triv}
\end{figure}

\begin{remark}
We are able to apply this construction for a non-reduced tree diagram.
However, the associated link diagram is different from the one of the equivalent reduced tree diagram:
Let $(T'_+, T'_-)$ be a non-reduced tree diagram obtained by inserting a caret into the reduced tree diagram $(T_+, T_-)$.
Then $\cL(T'_+, T'_-) = \cL(T_+, T_-) \sqcup \lower-1pt\hbox{$\bigcirc$}$ holds; see Figure \ref{caret_triv}.
For an element $g \in F$, let $\cL(g)$ denote the link diagram obtained from its reduced tree diagram.
\end{remark}

\subsection{3-colorablity of links}
In this section, we explain the 3-colorability of links which is one of fundamental invariants of links.

\begin{figure}[tbp]
\begin{center}
\includegraphics[height=60pt]{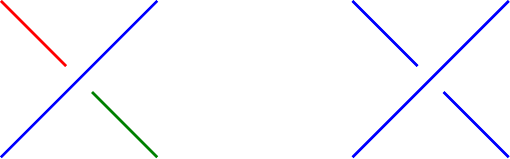}
\end{center}
\caption{The crossing with three colors and one color.}
\label{3-coloring}
\end{figure}

\begin{definition}
A \textbf{3-coloring} of a link diagram $D$ is a coloring with one of three colors (blue, green, and red) to each arc of $D$ such that the following condition is satisfied:
\begin{itemize}
\item at each crossing, three arcs are colored by either all different colors or the same color; see Figure \ref{3-coloring}. \label{3-color1}
\end{itemize}
A link diagram is \textbf{3-colorable} if there exists a 3-coloring with at least two colors.
\end{definition}

It is well known that the 3-colorability is a link invariant.

\begin{example}
The trefoil knot $3_1$ is 3-colorable.
The trivial knot $0_1$ and figure-eight knot $4_1$ are not 3-colorable.
\end{example}

\begin{figure}[tbp]
\begin{center}
\includegraphics[height=65pt]{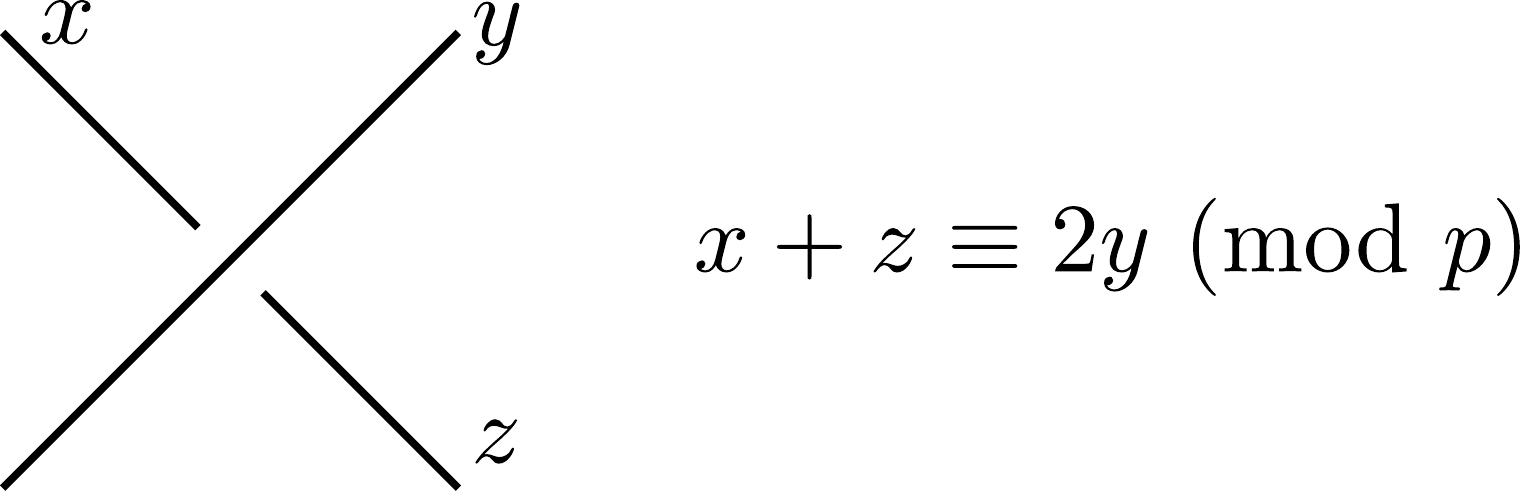}
\end{center}
\caption{The condition of the $p$-coloring.}
\label{p-coloring}
\end{figure}

\begin{remark} \label{p-col}
There exists a generalization of the 3-colorability as follows:
let $p$ be an integer grater than $2$.
A \textbf{$p$-coloring} of a link diagram $D$ is an assigning one of numbers $0, 1, \ldots, p-1$ to each arc of $D$ satisfying the following condition:
\begin{itemize}
\item at each crossing, the numbers $x$ and $z$ on two under-arcs, and number $y$ on one over-arc satisfy $x+z \equiv 2y \pmod p$; see Figure \ref{p-coloring}.
\end{itemize}
The \textbf{$p$-colorablility} of a link diagram is defined similarly.
The \textbf{coloring number} $c(L)$ of a link $L$ is the smallest $p$ such that $L$ is $p$-colorable.
\end{remark}

\section{Proof of Theorem \ref{main}}

We introduce a way of a 3-coloring of a link diagram $\cL(T_+, T_-)$ for any non-trivial element $(T_+, T_-)$ of $\cF$ to prove the main theorem directly.

\begin{figure}[tbp]
\begin{center}
\includegraphics[height=130pt]{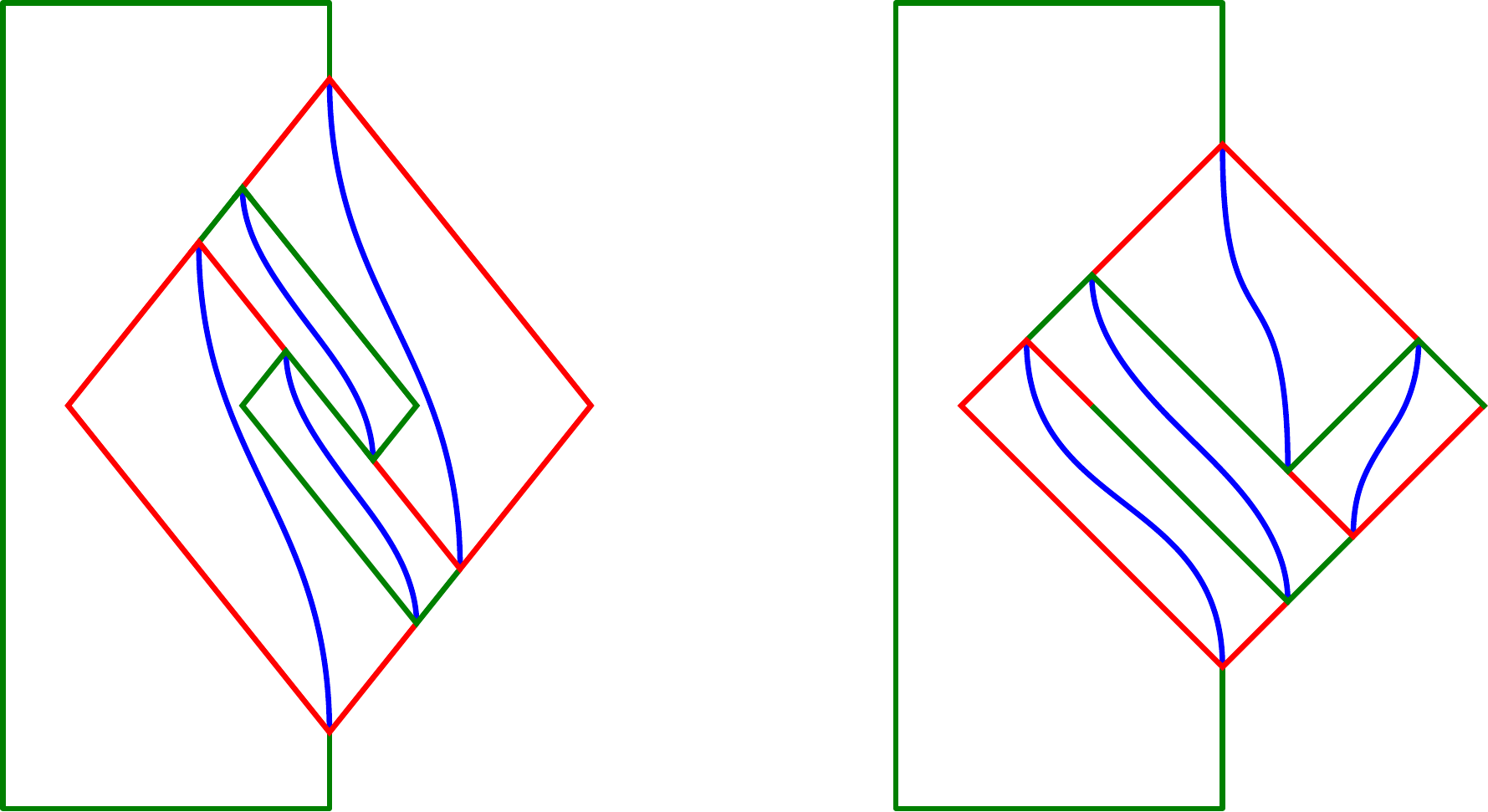}
\end{center}
\caption{Examples of the edge-coloring. The left figure induces a 3-coloring of the associated link diagram. However, the right one does not induce a 3-coloring since the second and fifth leaves have different colored edges.}
\label{edge-coloring}
\end{figure}

\begin{definition}
We color the edges of a graph $\cB(T_+, T_-)$ with three colors blue, green, and red as follows:
we color all the edges in the bounded regions of $\cB(T_+, T_-)$ with blue and the unique edge in the unbounded region with green.
The two edges around the root of $T_+$ are colored by two colors, thus we color the other edges with the distinguished color red.
The edges around the one of $T_-$ are colored similarly.
Then two edges around the children of the root are colored by two colors red and blue, thus we color the other edges with the distinguished color green.
The remaining edges of $\cB(T_+, T_-)$ are colored inductively; see Figure \ref{edge-coloring}.
We call this coloring the \textbf{edge-coloring} of $\cB(T_+, T_-)$.
\end{definition}

We induce a 3-coloring of $\cL(T_+, T_-)$ from the edge-coloring naturally.
By the rules of Jones' construction, it is clear that the induced arcs around each crossing are colored by three colors.
However, there generally exist two edges around a leaf such that their colors do not induce the coloring of an arc.
Therefore, it is sufficient to show that it does not occur for any element in $\cF$.

\begin{definition}
Let $(T_+, T_-)$ be in $\cT$.
We label the leaves of the trees with $1, \ldots, n$ as Section \ref{thompson}.
Since there exists the unique path from each leaf $i$ to the root of $T_+$ (resp.~ $T_-$), we define $\ell_i(T_+)$ (resp.~ $\ell_i(T_-)$) to be its length.
\end{definition}

\begin{remark}
In the edge-coloring of $\cB(T_+, T_-)$, the color of the edge in $T_+$ (resp.~ $T_-$) connecting the leaf $i$ depends only on the parity of $\ell_i(T_+)$ (resp.~ $\ell_i(T_-)$).
\end{remark}

From Theorem \ref{Thm_cF_Brown-Thompson}, the following holds:

\begin{lemma}[{\cite[Lemma 2.6]{aiello2021maximal}}]
Let $(T_+, T_-)$ be in $\cF$.
Then for each leaf $i$, $\ell_i(T_+) \equiv \ell_i(T_-) \pmod 2$ holds.
\end{lemma}

\begin{figure}[tbp]
\begin{center}
\includegraphics[height=150pt]{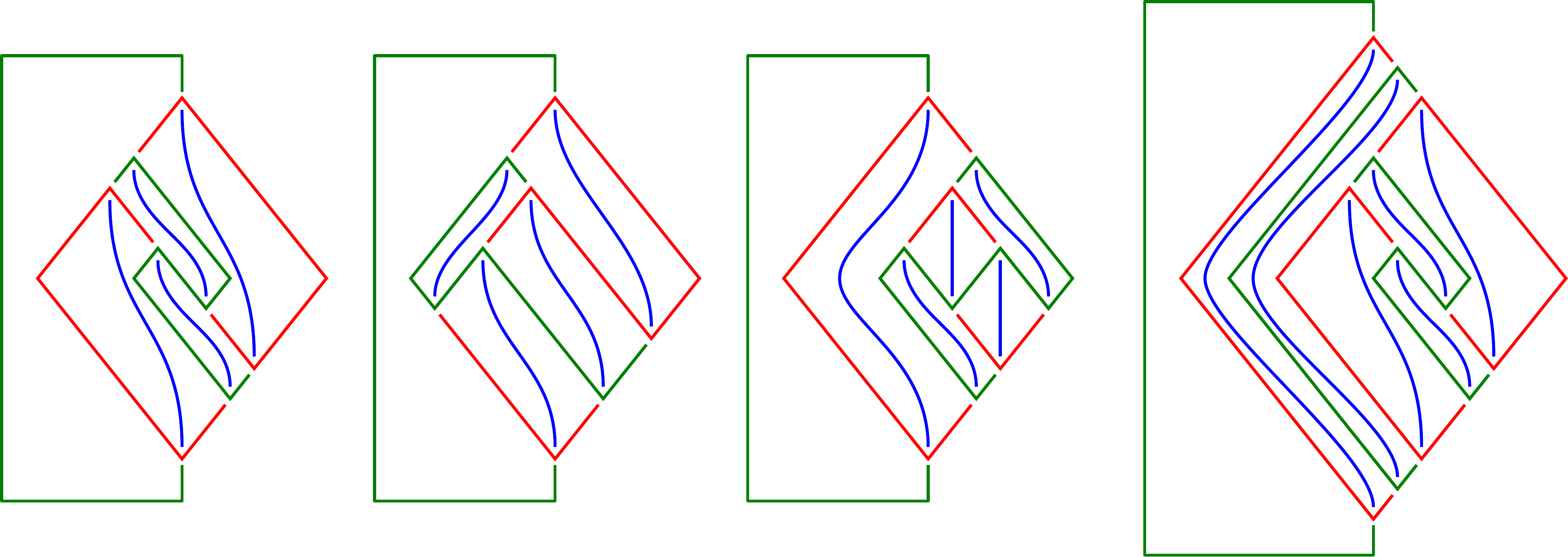}
\end{center}
\caption{Coloring of link diagrams of generators. }
\label{colF_3coloring}
\end{figure}

By above lemma, the edge-coloring of $\cB(T_+, T_-)$ in $\cF$ induces a 3-coloring of $\cL(T_+, T_-)$ with three colors.
Therefore Theorem \ref{main} immediately follows.
Figure \ref{colF_3coloring} gives the 3-coloring of the links associated with the generators of $\cF$.

Finally, there are natural questions about relationships between the properties of links and the 3-colorable subgroup $\cF$.

\begin{question}
Do the non-trivial elements of $\cF$ produce all 3-colorable knots and links?
\end{question}

Also we may consider a generalization of the 3-colorable subgroup:

\begin{question}
Does there exist a subgroup of $F$ whose non-trivial elements give $p$-colorable knots and links?; see Remark \ref{p-col}.
\end{question}

\section*{Acknowledgements}
We would like to thank Professor Tomohiro Fukaya who is the first author's supervisor for his several comments. 
We also wish to thank Professor Takuya Sakasai who is the second author's supervisor for his helpful comments.
We are also grateful to Professor Motoko Kato for her valuable suggestions.
We also are deeply grateful to Dr.\ Valeriano Aiello for his helpful advice.

\bibliographystyle{plain}
\bibliography{bib1} 
\bigskip
\address{
DEPARTMENT OF MATHEMATICAL SCIENCES,
TOKYO METROPOLITAN UNIVERSITY,
MINAMI-OSAWA HACHIOJI, TOKYO, 192-0397, JAPAN
}

\textit{E-mail address}: \href{mailto:kodama-yuya@ed.tmu.ac.jp}{\texttt{kodama-yuya@ed.tmu.ac.jp}}

\address{GRADUATE SCHOOL OF MATHEMATICAL SCIENCES, THE
UNIVERSITY OF TOKYO, 3-8-1 KOMABA, MEGURO-KU, TOKYO, 153-8914,
JAPAN}

\textit{E-mail address}: \href{mailto:takano@ms.u-tokyo.ac.jp}{\texttt{takano@ms.u-tokyo.ac.jp}}
\end{document}